\documentclass{article}

\usepackage{amsmath}
\usepackage{amssymb}
\usepackage{amsthm}

\numberwithin{equation}{section}
\newtheorem{thm}{Theorem}[section]

\newtheorem{prop}[thm]{Proposition}

\begin{document}
\title{
An inverse boundary value problem for the magnetic Schr\"{o}dinger operator with a bounded magnetic potential in a slab}
\author{Shitao Liu, Yang Yang}
\date{}
\maketitle

\begin{abstract}
We study an inverse boundary value problem with partial data in an infinite slab in $\mathbb{R}^{n}$, $n\geq 3$, for the magnetic Schr\"{o}dinger operator with an $L^{\infty}$ magnetic potential and an $L^{\infty}$ electric potential. We show that the magnetic field and the electric potential can be uniquely determined, when the Dirichlet and Neumann data are given on either different boundary hyperplanes or on the same boundary hyperplanes of the slab. This generalizes the result in \cite{KLU}, where the same uniqueness result was established when the magnetic potential is Lipschitz continuous. The proof is based on the complex geometric optics solutions constructed in \cite{KU}, which are special solutions to the magnetic Schr\"{o}dinger equation with $L^{\infty}$ magnetic and electric potentials in a bounded domain.
\end{abstract}

\section{INTRODUCTION AND STATEMENT OF RESULTS}
In this paper we study an inverse boundary value problems with partial data for the magnetic Schr\"{o}dinger operator with bounded magnetic potential and bounded electric potential in an infinite slab.

A magnetic Schr\"{o}dinger operator with bounded magnetic potential and bounded electric potential is a differential operator of the form
$$\mathcal{L}_{A,q}(x,D):=\displaystyle\sum^{n}_{j=1}(D_{j}+A_{j}(x))^{2}+q(x),$$
where $D_j=i^{-1}\frac{\partial}{\partial x_{j}}$, the complex-valued vector field $A=(A_{j})_{1\leq j\leq n}\in L^{\infty}(\Sigma;\mathbb{C}^{n})$ is the magnetic potential, and the complex-valued function $q\in L^{\infty}(\Sigma;\mathbb{C})$ is the electric potential. Throughout this paper, we shall assume $A$ and $q$ are compactly supported in the slab.

We denote the infinite slab by $\Sigma$. This is a subset of $\mathbb{R}^{n}$, $n\geq 3$, which lies between two parallel hyperplanes. By choosing appropriate coordinates, we may assume that the slab is defined by
$$\Sigma:=\{x=(x',x_{n})\in\mathbb{R}^{n}:x'=(x_{1},\dots,x_{n-1})\in\mathbb{R}^{n-1}, 0<x_{n}<L\},\quad L>0.$$
Its boundary consists of two parallel hyperplanes
$$\Gamma_{1}:=\{x\in\mathbb{R}^{n}:x_{n}=L\} \quad\quad \Gamma_{2}:=\{x\in\mathbb{R}^{n}:x_{n}=0\}.$$

Given $f\in H^{\frac{3}{2}}(\Gamma_{1})$ which is compactly supported in $\Gamma_{1}$, we are interested in the following Dirichlet problem
\begin{equation}\label{Dirichlet1}
\left\{
\begin{array}{rll}\vspace{1ex}
(\mathcal{L}_{A,q}-k^{2})u=0 \quad & \quad\textrm{ in } \Sigma \\ \vspace{1ex}
u=f \quad & \quad\textrm{ on } \Gamma_{1} \\ \vspace{1ex}
u=0 \quad & \quad\textrm{ on } \Gamma_{2}. \\
\end{array}
\right.
\end{equation}
Here $k\geq 0$ is fixed. This Dirichlet boundary value problem admits a unique solution in $H^{2}_{loc}(\overline{\Sigma})$ for admissible $k$; in this case, we define the Dirichlet-to-Neumann map by
$$
\begin{array}{rcl}\vspace{1ex}
\mathcal{N}_{A,q}: \quad (H^{\frac{3}{2}}(\Gamma_{1})\cap\mathcal{E}'(\Gamma_{1})) & \rightarrow & H^{\frac{1}{2}}_{loc}(\partial\Sigma)\\ \vspace{1ex}
f & \mapsto & (\partial_{\nu}+iA\cdot\nu)u|_{\partial\Sigma},
\end{array}
$$
where $u$ is the solution of \eqref{Dirichlet1}, $\mathcal{E}'(\Gamma_{1})$ is the set of compactly supported distributions on $\Gamma_{1}$, $\nu$ is the unit outer normal vector field to $\partial\Sigma=\Gamma_{1}\cup\Gamma_{2}$.

It was observed in \cite{S} that the Dirichlet-to-Neumann map is invariant under a certain gauge transformation of the magnetic potential. More precisely, from the identities
\begin{equation}\label{gauge}
e^{-i\Psi}\mathcal{L}_{A,q}e^{i\Psi}=\mathcal{L}_{A+\nabla\Psi,q},\quad e^{-i\Psi}\mathcal{N}_{A,q}e^{i\Psi}=\mathcal{N}_{A+\nabla\Psi,q}
\end{equation}
we conclude that $\mathcal{N}_{A,q}=\mathcal{N}_{A+\nabla\Psi,q}$ provided that $\Psi\in W^{1,\infty}(\overline{\Sigma})$ is compactly supported and $\Psi|_{\partial\Sigma}=0$. Thus, $\mathcal{N}_{A,q}$ only gives information about the magnetic field $dA$, where the vector field $A$ is viewed as the 1-form $\sum^{n}_{j=1}A_{j}dx_{j}$, and the differential 2-form $dA$ is defined as
$$dA=\displaystyle\sum_{1\leq j<k\leq n}(\partial_{x_{j}}A_{k}-\partial_{x_{k}}A_{j})dx_{j}\wedge dx_{k}$$
in the sense of distributions.

In this paper we consider the inverse boundary value problem of determining the magnetic field $dA$ and the electric potential $q$ for only bounded $A$ and $q$. Our main results in this paper are the following two theorems. They generalize the results in \cite{LU}, which were obtained when $A=0$, and the results in \cite{KLU}, which were obtained when $A\in W^{1,\infty}(\Sigma;\mathbb{C}^{n})$ is compactly supported. The first result concerns the case when the data and the measurements are on different boundary hyperplanes, while the second deals with the case when the data and the measurements are on the same boundary hyperplanes.

\begin{thm}
Let $\Sigma\subset\mathbb{R}^{n}(n\geq 3)$ be an infinite slab with boundary hyperplanes $\Gamma_{1}$ and $\Gamma_{2}$. Let $A^{(j)}\in L^{\infty}({\Sigma};\mathbb{C}^{n})\cap\mathcal{E}'(\bar{\Sigma};\mathbb{C}^{n})$ and $q^{(j)}\in L^{\infty}(\Sigma;\mathbb{C})\cap\mathcal{E}'(\bar{\Sigma};\mathbb{C}), j=1,2.$ Denote by $B\subset\mathbb{R}^{n}$ an open ball containing the supports of $A^{(j)},q^{(j)},j=1,2$. Denote by $\gamma_{j}\subset\Gamma_{j}$ an arbitrary open subset of $\Gamma_{j}$ satisfying that
$$\Gamma_{j}\cap\bar{B}\subset\gamma_{j},\quad j=1,2.$$
Let $k\geq 0$ be fixed and admissible for both the operator $\mathcal{L}_{A^{(j)},q^{(j)}}$ and its real transpose $\mathcal{L}_{-A^{(j)},q^{(j)}}$, $j=1,2$. If
$$\mathcal{N}_{A^{(1)},q^{(1)}}(f)|_{\gamma_{2}}=\mathcal{N}_{A^{(2)},q^{(2)}}(f)|_{\gamma_{2}}$$
for all $f\in H^{\frac{3}{2}}(\Gamma_{1})$ with $supp(f)\subset\gamma_{1}$, then $dA^{(1)}=dA^{(2)}$ and $q^{(1)}=q^{(2)}$.
\end{thm}

The assumption that $k\geq 0$ is admissible for the real transpose $\mathcal{L}_{-A^{(j)},q^{(j)}}$ is necessary for the proof of a Runge type approximation result. A special case of the assumption in the above theorem is that the supports of $A^{(j)},q^{(j)}$ are strictly contained in the interior of the slab, in this situation $B\cap\Gamma_1 = B\cap\Gamma_2 = \emptyset$, so the open subsets $\gamma_{1}$ and $\gamma_{2}$ can be arbitrarily small in $\Gamma_1$ and $\Gamma_2$.

Next we give the result of our inverse problem for the case when the given data and the measurements are on the same hyperplane $\Gamma_1$.

\begin{thm}
Let $\Sigma\subset\mathbb{R}^{n}(n\geq 3)$ be an infinite slab with boundary hyperplanes $\Gamma_{1}$ and $\Gamma_{2}$. Let $A^{(j)}\in L^{\infty}({\Sigma};\mathbb{C}^{n})\cap\mathcal{E}'(\bar{\Sigma};\mathbb{C}^{n})$ and $q^{(j)}\in L^{\infty}(\Sigma;\mathbb{C})\cap\mathcal{E}'(\bar{\Sigma};\mathbb{C}), j=1,2.$ Denote by $B\subset\mathbb{R}^{n}$ an open ball containing the supports of $A^{(j)},q^{(j)},j=1,2$. Denote by $\gamma_{1},\gamma'_{1}\subset\Gamma_{1}$ two arbitrary open subsets of $\Gamma_{1}$ satisfying that
$$\Gamma_{1}\cap\bar{B}\subset\gamma_{1},\quad \Gamma_{1}\cap\bar{B}\subset\gamma'_{1}.$$
Let $k\geq 0$ be fixed and admissible for both the operator $\mathcal{L}_{A^{(j)},q^{(j)}}$ and its real transpose $\mathcal{L}_{-A^{(j)},q^{(j)}}$, $j=1,2$. If
$$\mathcal{N}_{A^{(1)},q^{(1)}}(f)|_{\gamma'_{1}}=\mathcal{N}_{A^{(2)},q^{(2)}}(f)|_{\gamma'_{1}}$$
for all $f\in H^{\frac{3}{2}}(\Gamma_{1})$ with $supp(f)\subset\gamma_{1}$, then $dA^{(1)}=dA^{(2)}$ and $q^{(1)}=q^{(2)}$.
\end{thm}

The proofs of Theorem 1.1 and Theorem 1.2 are based on the construction of a special kind of complex geometric optics (CGO) solutions which vanish on appropriate subsets of the boundary hyperplanes, through a reflection argument. The constructions of CGO solutions for the conductivity equation and the Schr\"{o}dinger equation go back to \cite{SU}. Constructing CGO solutions by using a reflection argument was initiated in \cite{I}. It was applied to the inverse boundary value problem for the Schr\"{o}dinger operator in a slab in \cite{LU}, and for the magnetic Sch\"{o}dinger operator in a slab in \cite{KLU}.

The geometry of an infinite slab shows up in modeling the propagation of acoustic waves in shallow oceans, see \cite{AK}. Various inverse problems in a slab have been studied by many authors. In \cite{Ik,SW}, inverse conductivity problems in a slab were studied; in \cite{CGI}, an inverse scattering problem for the Schr\"{o}dinger operator in a slab was considered.

Another closely related inverse problem for the magnetic Schr\"{o}dinger operator is the inverse boundary value problem in a bounded open subset $\Omega$ of $\mathbb{R}^{n}$, $n\geq 3$. In this problem one hopes to determine the magnetic field and the electric potential by measuring the Dirichlet-to-Neumann map on the boundary $\partial\Omega$, under various regularity assumption on the magnetic potential $A$ and the electric potential $q$. In \cite{S}, this type of determination was established for $W^{2,\infty}$ magnetic potentials, which are small in a suitable sense, and $L^{\infty}$ electric potentials. In \cite{NSU}, the smallness condition was removed for $C^{\infty}$ magnetic and electric potentials. This uniqueness was extended to $C^{1}$ magnetic potentials in \cite{T2}, to some less regularity but small potentials in \cite{P}, to Dini continuous magnetic potentials in \cite{S1}, and recently to $L^{\infty}$ magnetic and electric potentials in \cite{KU}.

In a bounded domain $\Omega$, one may also consider this inverse boundary value problem with partial boundary measurements. More precisely, let $\gamma_{1}$ and $\gamma_{2}$ be two open subsets of $\partial\Omega$, and let the partial data measurements be $\mathcal{N}_{A,q}(f)|_{\gamma_{2}}$ for all $f$ with $supp(f)\subset\gamma_{1}$. Determination of the magnetic potential and the electric potential from partial boundary measurements was proved in \cite{FKSU} when $\gamma_{1}=\partial\Omega$ and $\gamma_{2}$ is possibly a very small subset. Under the assumption that $A^{(1)}=A^{(2)}$ and $q^{(1)}=q^{(2)}$ in a neighborhood of the boundary, it was proved in \cite{BJ} when $\gamma_{1}=\partial\Omega$ and $\gamma_{2}$ is arbitrary. This was generalized to arbitrary non-empty subsets $\gamma_{1},\gamma_{2}\subset\partial\Omega$ in \cite{KLU}. Logarithmic stability estimate for this problem was obtained in \cite{BJ}.

This paper is structured as follows. In Section 2, we review the construction of CGO solutions in a bounded domain $\Omega\subset\mathbb{R}^{n}$, $n\geq 3$, to the magnetic Schr\"{o}dinger equation with $L^{\infty}$ magnetic and electric potentials, following \cite{KU}. Section 3 contains the proof of Theorem 1.1, and is divided into three parts: Section 3.1 establishes an integral inequality of Green's type for some special solutions of the magnetic Schr\"{o}dinger equation; Section 3.2 shows how to construct the desired partial data CGO solutions from the CGO solutions constructed in Section 2 through a reflection argument; In Section 3.3, we finish the proof of Theorem~1.1 by combining the result in Section 3.1 and 3.2. Section 4 demonstrates the proof of Theorem 1.2. 

\section{CGO SOLUTIONS}
In this section, we recall how to construct CGO solutions for the mangetic Schr\"{o}dinger operator with a bounded magnetic potential on a bounded domain. For detailed construction we refer to \cite{KU}. This construction together with a reflection argument will allow us to construct the special partial data CGO solutions in the infinite slab $\Sigma$.

Let $\Omega\subset\mathbb{R}^{n}$, $n\geq 3$, be a bounded open subset with smooth boundary. Rewrite the magnetic Schr\"{o}dinger equation as
$$
\begin{array}{rl}
\mathcal{L}_{A,q}u & =-\Delta u(x)+A(x)\cdot Du(x)+D\cdot(A(x)u(x))+(A^{2}(x)+q(x))u(x)\\
 &=0 \quad\quad\quad\quad \textrm{ in }\Omega,
\end{array}
$$
with $A\in L^{\infty}(\Omega;\mathbb{C}^{n})$ and $q\in L^{\infty}(\Omega;\mathbb{C})$. Notice that $Au\in L^{\infty}(\Omega;\mathbb{C}^{n})\cap\mathcal{E}'(\Omega;\mathbb{C}^{n})$ and $D\cdot(A(x)u(x))\in H^{-1}(\Omega)$, therefore
$$\mathcal{L}_{A,q}: C^{\infty}_{c}(\Omega)\rightarrow H^{-1}(\mathbb{R}^{n})\cap\mathcal{E}'(\Omega)$$
is a bounded operator. Here $\mathcal{E}'(\Omega):=\{u\in\mathcal{D}'(\Omega): supp(u) \textrm{ is compact}\}$. Following \cite{KU}, we will construct solutions to the above magnetic Schr\"{o}dinger equation of the form
\begin{equation}\label{CGO}
u(x,\zeta,h)=e^{x\cdot\zeta/h}(a(x,\zeta,h)+r(x,\zeta,h)).
\end{equation}
Here $\zeta\in\mathbb{C}^{n}$ is a complex vector with $\zeta\cdot\zeta=0$. In this paper we shall work with $\zeta$ which can be written as $\zeta=\zeta^{(0)}+\zeta^{(1)}$ with $\zeta^{(0)}$ independent of $h$, $\zeta^{(0)}\cdot\zeta^{(0)}=0$, $|Re\zeta^{(0)}|=|Im\zeta^{(0)}|=1$, and $\zeta^{(1)}=\mathcal{O}(h)$ as $h\rightarrow 0$. $a$ is a smooth function, $r$ is an $L^{2}(\Omega)$-function, and $h>0$ is a small semiclassical parameter. They are called complex geometric optics solutions since the phase functions are complex-valued. This construction relies on a Carleman estimate for the magnetic Schr\"{o}dinger operator $\mathcal{L}_{A,q}$ with a gain of two derivatives, which is based on the corresponding Carleman estimate for the Laplacian obtained in \cite{ST}.

We extend $A$ to a vector field on $\mathbb{R}^{n}$ by defining it to be zero in $\mathbb{R}^{n}\backslash\Omega$, this extension is still denoted by $A$. Then $A\in (L^{\infty}\cap\mathcal{E}')(\mathbb{R}^{n};\mathbb{C}^{n})\subset L^{p}(\mathbb{R}^{n};\mathbb{C}^{n})$, $1\leq p\leq\infty$. Let $\eta_{\epsilon}(x)=\epsilon^{-n}\eta(x/\epsilon),\epsilon>0$ be the standard mollifier with $\eta\in C^{\infty}_{c}(\mathbb{R}^{n})$, $0\leq\eta\leq 1$ and $\int\eta\,dx=1$. Define $A^{\sharp}:=A\ast\eta_{\epsilon}\in C^{\infty}_{c}(\mathbb{R}^{n};\mathbb{C}^{n})$, then
\begin{equation}\label{Asharp}
\|A-A^{\sharp}\|_{L^{2}(\mathbb{R}^{n})}=\mathcal{O}(\epsilon),\quad \textrm{ as }\epsilon\rightarrow 0.
\end{equation}
and
\begin{equation}\label{Asharpderivative}
\|\partial^{\alpha}A^{\sharp}\|_{L^{\infty}(\mathbb{R}^{n})}=\mathcal{O}(\epsilon^{-|\alpha|}), \quad\textrm{ as }\epsilon\rightarrow 0 \textrm{ for all } \alpha, |\alpha|\geq 0.
\end{equation}

Consider the conjugated operator
\begin{equation}\label{conjugate}
\begin{array}{rl}
e^{-x\cdot\zeta/h}h^{2}\mathcal{L}_{A,q}e^{x\cdot\zeta/h}=-h^{2}\Delta-2ih\zeta^{(0)}\cdot D-2ih\zeta^{(1)}\cdot D+h^{2}A\cdot D-2ih\zeta^{(0)}\cdot A^{\sharp}\\
-2ih\zeta^{(0)}\cdot (A-A^{\sharp})-2ih\zeta^{(1)}\cdot A+h^{2}D\cdot(A\cdot)+h^{2}(A^{2}+q).\\
\end{array}
\end{equation}
We choose $a$ so that it satisfies the equation
$$
\zeta^{(0)}\cdot Da+\zeta^{(0)}\cdot A^{\sharp}a=0 \quad\quad \textrm{ in } \mathbb{R}^{n}.
$$
This is a transport equation. To solve this equation we look for solutions of the form $a=e^{\Phi^{\sharp}}$ where $\Phi^{\sharp}$ solves the equation
\begin{equation}\label{Phi}
\zeta^{(0)}\cdot D\Phi^{\sharp}+\zeta^{(0)}\cdot A^{\sharp}=0 \quad\quad \textrm{ in } \mathbb{R}^{n}.
\end{equation}
As $\zeta^{(0)}\cdot\zeta^{(0)}=0$ and $|Re\,\zeta^{(0)}|$=$|Im\,\zeta^{(0)}|=1$, the operator $N_{\zeta^{(0)}}:=\zeta^{(0)}\cdot\nabla$ is the $\bar{\partial}$-operator in appropriate coordinates. To solve the equation \eqref{Phi}, we introduce the Cauchy transform $N^{-1}_{\zeta^{(0)}}$:
$$(N^{-1}_{\zeta^{(0)}}f)(x):=\displaystyle\frac{1}{2\pi}\displaystyle\int_{\mathbb{R}^{n}}\frac{f(x-y_{1} Re\,\zeta^{(0)}-
y_{2} Im\,\zeta^{(0)})}{y_{1}+iy_{2}}\,dy_{1}dy_{2}, \quad f\in C_{c}(\mathbb{R}^{n}).$$
By \cite[Lemma 2.4, 2.5]{KU}, $\Phi^{\sharp}=N^{-1}_{\zeta^{(0)}}(-i\zeta^{(0)}\cdot A^{\sharp})\in C^{\infty}(\mathbb{R}^{n})$ solves the equation \eqref{Phi}, and satisfies
\begin{equation}\label{Phiderivative}
\|\partial^{\alpha}\Phi^{\sharp}\|_{L^{\infty}(\mathbb{R}^{n})}=\mathcal{O}(\epsilon^{-|\alpha|}), \quad\textrm{ as }\epsilon\rightarrow 0 \textrm{ for all } \alpha, |\alpha|\geq 0.
\end{equation}
Furthermore, if we denote $\Phi(\cdot,\zeta^{(0)}):=N^{-1}_{\zeta^{(0)}}(-i\zeta^{(0)}\cdot A)\in L^{\infty}(\mathbb{R}^{n})$, then
\begin{equation}\label{convergence}
\Phi^{\sharp}(\cdot,\zeta^{(0)},\epsilon)\rightarrow\Phi(\cdot,\zeta^{(0)}) \textrm{ in } L^{2}_{loc}(\mathbb{R}^{n}) \textrm{ as } \epsilon\rightarrow 0.
\end{equation}

To make \eqref{CGO} a solution of the equation $\mathcal{L}_{A,q}u=0$, the correction term $r$ needs to satisfy the equation
$$e^{-x\cdot\zeta/h}h^{2}\mathcal{L}_{A,q}e^{x\cdot\zeta/h}r=-e^{-x\cdot\zeta/h}h^{2}\mathcal{L}_{A,q}e^{x\cdot\zeta/h}a.$$
Notice that the right hand side belongs to $H^{-1}(\Omega)$, and is of order $o(h)$ due to our choice of $a$. Hence $r$ can be found using the solvability result \cite[Proposition 2.3]{KU}; moreover, $r$ satisfies the decay property $\|r\|_{H^{1}_{scl}(\Omega)}=\mathcal{O}(h\epsilon^{-2}+\epsilon)$ as $h\rightarrow 0$, where $\|r\|^{2}_{H^{1}_{scl}(\Omega)}=\|r\|^{2}_{L^{2}(\Omega)}+\|hD r\|^{2}_{L^{2}(\Omega)}$. Putting these together and choosing $\epsilon=h^{1/3}$, we obtain the following result (see \cite[Proposition 2.6]{KU})
\begin{prop}\label{existence}
Let $\Omega\subset\mathbb{R}^{n}$, $n\geq 3$, be a bounded open set. Let $A\in L^{\infty}(\Omega;\mathbb{C}^{n})$ and $q\in L^{\infty}(\Omega;\mathbb{C})$. Let $\zeta\in\mathbb{C}^{n}$ be such that $\zeta\cdot\zeta=0$, $\zeta=\zeta^{(0)}+\zeta^{(1)}$ with $\zeta^{(0)}$ independent of $h$, $|Re\,\zeta^{(0)}|=|Im\,\zeta^{(0)}|=1$, and $\zeta^{(1)}=\mathcal{O}(h)$ as $h\rightarrow 0$. Then for $h>0$ small, there exist solutions $u(x,\zeta,h)\in H^{1}(\Omega)$ to the magnetic Schr\"{o}dinger equation $\mathcal{L}_{A,q}u=0$ of the form
$$u(x,\zeta,h)=e^{x\cdot\zeta/h}(e^{\Phi^{\sharp}(x,\zeta^{(0)},h)}+r(x,\zeta,h)).$$
Here the function $\Phi^{\sharp}(\cdot,\zeta^{(0)},h)\in C^{\infty}(\mathbb{R}^{n})$ satisfies $\|\partial^{\alpha}\Phi^{\sharp}\|_{L^{\infty}(\mathbb{R}^{n})}\leq C_{\alpha}h^{-|\alpha|/3}$ for all $\alpha$ with $|\alpha|\geq 0$, and $\Phi^{\sharp}(\cdot,\zeta^{(0)},h)$ converges to $\Phi(\cdot,\zeta^{(0)}):=N^{-1}_{\zeta^{(0)}}(-i\zeta^{(0)}\cdot A)\in L^{\infty}(\mathbb{R}^{n})$ in $L^{2}_{loc}(\mathbb{R}^{n})$ as $h\rightarrow 0$. Here we have extended $A$ by zero to $\mathbb{R}^{n}\backslash\Omega$. The remainder $r$ satisfies $\|r\|_{H^{1}_{scl}(\Omega)}=\mathcal{O}(h^{1/3})$ as $h\rightarrow 0$.
\end{prop}

In the proofs below we may need CGO solutions belonging to $H^{2}(\Omega)$. To obtain such solutions, we can first construct CGO solutions on a larger bounded open set whose interior contains the closure of $\Omega$, and then restrict the solutions back to $\Omega$. These solutions will belong to $H^{2}(\Omega)$ by elliptic regularity.

\section{PROOF OF THEOREM 1.1}

\subsection{INTEGRAL IDENTITY}
First we derive an integral identity from Green's formula. We shall use the $L^{2}$-spaces with inner products
$$(u,v)_{L^{2}(\Omega)}=\displaystyle\int_{\Omega}u(x)\overline{v(x)}\,dx \quad (u,v)_{L^{2}(\partial\Omega)}=\displaystyle\int_{\partial\Omega}u(x)\overline{v(x)}\,dS $$
where $dS$ is the surface measure on $\partial\Omega$. The following Green's formula for the magnetic Schr\"{o}dinger operator $\mathcal{L}_{A,q}$ was established in \cite{FKSU}.
\begin{equation}\label{green}
(\mathcal{L}_{A,q}u,v)_{L^{2}(\Omega)}-(u,\mathcal{L}_{\bar{A},\bar{q}}v)_{L^{2}(\Omega)}=
(u,(\partial_{\nu}+i\nu\cdot\bar{A})v)_{L^{2}(\partial\Omega)}-((\partial_{\nu}+i\nu\cdot A)u,v)_{L^{2}(\partial\Omega)}
\end{equation}
for all $u,v\in H^{2}(\Omega)$.

Let $k\geq 0$ be fixed and admissible for both the operator $\mathcal{L}_{A^{(j)},q^{(j)}}$ and its real transpose $\mathcal{L}_{-A^{(j)},q^{(j)}}$, $j=1,2$. Let $u_{1}\in H^{2}_{loc}(\overline{\Sigma})$ be the admissible solution to the Dirichlet boundary value problem
$$
\left\{
\begin{array}{rll}\vspace{1ex}
(\mathcal{L}_{A^{(1)},q^{(1)}}-k^{2})u_{1}=0 \quad & \quad\textrm{ in } \Sigma \\ \vspace{1ex}
u_{1}=f \quad & \quad\textrm{ on } \Gamma_{1} \\ \vspace{1ex}
u_{1}=0 \quad & \quad\textrm{ on } \Gamma_{2} \\
\end{array}
\right.
$$
where $f\in H^{3/2}(\Gamma_{1})$ with $supp(f)\subset\gamma_{1}$. Let $w\in H^{2}_{loc}(\overline{\Sigma})$ be the admissible solution of the problem
$$
\left\{
\begin{array}{rll}\vspace{1ex}
(\mathcal{L}_{A^{(1)},q^{(1)}}-k^{2})w=0 \quad & \quad\textrm{ in } \Sigma \\ \vspace{1ex}
w=u_{1} \quad & \quad\textrm{ on } \Gamma_{1} \\ \vspace{1ex}
w=u_{1} \quad & \quad\textrm{ on } \Gamma_{2}. \\
\end{array}
\right.
$$
Set $v=w-u_{1}$, then
\begin{equation}\label{difference}
\begin{array}{rl}
(\mathcal{L}_{A^{(2)},q^{(2)}}-k^{2})v & =(A^{(1)}-A^{(2)})\cdot Du_{1}+D\cdot((A^{(1)}-A^{(2)})u_{1})\\
 & +((A^{(1)})^{2}-(A^{(2)})^{2}+q^{(1)}-q^{(2)})u_{1} \quad \textrm{ in } \Sigma.
\end{array}
\end{equation}
Under the assumption that $\mathcal{N}_{A^{(1)},q^{(1)}}(f)|_{\gamma_{2}}=\mathcal{N}_{A^{(2)},q^{(2)}}(f)|_{\gamma_{2}}$, we have
$$(\partial_{\nu}+iA^{(1)}\cdot\nu)u_{1}|_{\gamma_{2}}=(\partial_{\nu}+iA^{(2)}\cdot\nu)w|_{\gamma_{2}}.$$
Since $u_{1}=w=0$ on $\Gamma_{2}$, we conclude that $(\partial_{\nu}v)|_{\gamma_{2}}=0$. Introduce the notations
$$l_{1}:=\Gamma_{1}\cap\bar{B}\subset\gamma_{1}, \quad l_{2}:=\Gamma_{2}\cap\bar{B}\subset\gamma_{2}, \quad l_{3}:=\Sigma\cap\partial B.$$
It is clear that $\partial(\Sigma\cap B)=l_{1}\cup l_{2}\cup l_{3}$. Then $v\in H^{2}_{loc}(\overline{\Sigma})$ is a solution to the equation
$$(-\Delta-k^{2})v=0 \quad\quad \textrm{ in } \Sigma\backslash\bar{B}$$
with $v=\partial_{\nu}v=0$ on $\gamma_{2}\backslash\overline{l_{2}}$, thus by unique continuation, $v=0$ in $\Sigma\backslash\bar{B}$. As a consequence, $v=\partial_{\nu}v=0$ on $l_{3}$.

Let $u_{2}\in H^{2}(\Sigma\cap B)$ be a solution to the equation
\begin{equation}\label{u2}
(\mathcal{L}_{\overline{A^{(2)}},\overline{q^{(2)}}}-k^{2})u_{2}=0 \quad\quad \textrm{ in } \Sigma\cap B
\end{equation}
with $u_{2}=0 \textrm{ on } l_{1}$. Apply Green's formula \eqref{green} to $v$ and $u_{2}$ over $\Sigma\cap B$ to get
\begin{equation}\label{green2}
\begin{array}{rl}
  & ((\mathcal{L}_{A^{(2)},q^{(2)}}-k^{2})v,u_{2})_{L^{2}(\Sigma\cap B)}-(v,(\mathcal{L}_{\overline{A^{(2)}},\overline{q^{(2)}}}-k^{2})u_{2})_{L^{2}(\Sigma\cap B)}\\
 =& (v,(\partial_{\nu}+i\nu\cdot\overline{A^{(2)}})u_{2})_{L^{2}(\partial(\Sigma\cap B))}-((\partial_{\nu}+i\nu\cdot A^{(2)})v,u_{2})_{L^{2}(\partial(\Sigma\cap B))}.
\end{array}
\end{equation}
Notice that $v=0$ on $l_{1}\cup l_{2}\cup l_{3}$, $\partial_{\nu}v=0$ on $l_{2}\cup l_{3}$ and $u_{2}=0$ on $l_{1}$, the equation \eqref{green2} then reduces to
$$((\mathcal{L}_{A^{(2)},q^{(2)}}-k^{2})v,u_{2})_{L^{2}(\Sigma\cap B)}=0.$$
This together with \eqref{difference} gives, after integrating by parts, that
\begin{equation}\label{identity1}
\begin{array}{rl}\vspace{1ex}
& \displaystyle\int_{\Sigma\cap B}(A^{(1)}-A^{(2)})\cdot ((Du_{1})\overline{u_{2}}+u_{1}\overline{Du_{2}})\,dx+\displaystyle\frac{1}{i}\displaystyle\int_{\partial(\Sigma\cap B)}(A^{(1)}-A^{(2)})\cdot\nu u_{1}\overline{u_{2}}\,dS\\ \vspace{1ex}
+& \displaystyle\int_{\Sigma\cap B}((A^{(1)})^{2}-(A^{(2)})^{2}+q^{(1)}-q^{(2)})u_{1}\overline{u_{2}}\,dx=0.\\
\end{array}
\end{equation}

We would like to show that the second term on the left vanishes. Without loss of generality, we may assume that
\begin{equation}\label{normal}
A^{(1)}\cdot\nu=A^{(2)}\cdot\nu=0 \quad \textrm{ on } \Gamma_{1}\cup\Gamma_{2}.
\end{equation}
In fact, if this does not hold, we can find compactly supported $\Psi^{(j)}\in W^{1,\infty}(\overline{\Sigma})$, $j=1,2$, such that
$$\Psi^{(j)}|_{\partial\Sigma}=0, \quad \textrm{ and } \partial_{\nu}\Psi^{(j)}=-A^{(j)}\cdot\nu \textrm{ on } \partial\Sigma.$$
Then we can replace $A^{(j)}$ by $A^{(j)}+\nabla\Psi^{(j)}$. This will not change the Dirichlet-to-Neumann map due to \eqref{gauge}. The existence of such $\Psi^{(j)}$ was proved in \cite[Theorem 1.3.3]{H}. On the other hand, we conclude $A^{(1)}=A^{(2)}=0$ on $l_{3}$, since their supports are contained in $B$. Putting these together, we have showed that $(A^{(1)}-A^{(2)})\cdot\nu=0$ on the boundary $\partial(\Sigma\cap B)$, thus the second term on the left hand side in \eqref{identity1} vanishes, and we obtain
\begin{equation}\label{identity2}
\begin{array}{rl}\vspace{1ex}
& \displaystyle\int_{\Sigma\cap B}(A^{(1)}-A^{(2)})\cdot ((Du_{1})\overline{u_{2}}+u_{1}\overline{Du_{2}})\,dx\\ \vspace{1ex}
+& \displaystyle\int_{\Sigma\cap B}((A^{(1)})^{2}-(A^{(2)})^{2}+q^{(1)}-q^{(2)})u_{1}\overline{u_{2}}\,dx=0.\\
\end{array}
\end{equation}

Introduce the following function spaces
$$
\begin{array}{rl} \vspace{1ex}
\mathcal{W}(\Sigma):=&\{u\in H^{2}_{loc}(\overline{\Sigma}) \textrm{ admissible}:(\mathcal{L}_{A^{(1)},q^{(1)}}-k^{2})u=0 \textrm{ in } \Sigma,\\ 
 & u|_{\Gamma_{2}}=0, supp(u|_{\Gamma_{1}})\subset\gamma_{1}, u \}\\
\end{array}
$$
$$
\mathcal{V}_{l_{j}}(\Sigma\cap B):=\{u\in H^{2}(\Sigma\cap B):(\mathcal{L}_{\overline{A^{(2)}},\overline{q^{(2)}}}-k^{2})u=0 \textrm{ in } \Sigma\cap B, u|_{l_{j}}=0\}.
$$
$$
\mathcal{W}_{l_{j}}(\Sigma\cap B):=\{u\in H^{2}(\Sigma\cap B):(\mathcal{L}_{A^{(1)},q^{(1)}}-k^{2})u=0 \textrm{ in } \Sigma\cap B, u|_{l_{j}}=0\}.
$$
Then the above argument shows that identity \eqref{identity2} holds for all $u_{1}\in \mathcal{W}(\Sigma)$ and for all $u_{2}\in \mathcal{V}_{l_{1}}(\Sigma\cap B)$. The following proposition allows us to enlarge the function space where $u_{1}$ lies from $\mathcal{W}(\Sigma)$ to $\mathcal{W}_{l_{2}}(\Sigma\cap B)$.
\begin{prop}
$\mathcal{W}(\Sigma)$ is a dense subspace of $\mathcal{W}_{l_{2}}(\Sigma\cap B)$ in the $L^{2}(\Sigma\cap B)$-topology.
\end{prop}
\noindent This is a Runge type approximation result. The proof is exactly the same as the one for \cite[Proposition 3.1]{KLU}, since in that proof no regularity assumption of $A^{(j)}$ and $q^{(j)}$ is involved.

Summing up, we have proved
\begin{prop}
With above notations, the identity
\begin{equation}\label{identity3}
\begin{array}{rl}\vspace{1ex}
& \displaystyle\int_{\Sigma\cap B}(A^{(1)}-A^{(2)})\cdot ((Du_{1})\overline{u_{2}}+u_{1}\overline{Du_{2}})\,dx\\ \vspace{1ex}
+& \displaystyle\int_{\Sigma\cap B}((A^{(1)})^{2}-(A^{(2)})^{2}+q^{(1)}-q^{(2)})u_{1}\overline{u_{2}}\,dx=0\\
\end{array}
\end{equation}
holds for all $u_{1}\in W_{l_{2}}(\Sigma\cap B)$ and for all $u_{2}\in V_{l_{1}}(\Sigma\cap B)$.
\end{prop}

\subsection{PARTIAL DATA CGO SOLUTIONS}
Next we construct partial data CGO solutions $u_{1}\in \mathcal{W}_{l_{2}}(\Sigma\cap B)$ and $u_{2}\in \mathcal{V}_{l_{1}}(\Sigma\cap B)$. The idea is to first construct CGO solutions on a larger domain containing $\Sigma\cap B$, and then take the difference of a solution and its reflection to fulfill the boundary condition. Recall that in the construction of CGO solutions on a bounded domain, we have used complex vectors $\zeta\in\mathbb{C}^{n}$ satisfying that $\zeta\cdot\zeta=0$, $\zeta=\zeta^{(0)}+\zeta^{(1)}$ with $\zeta^{(0)}$ independent of $h$, $|Re\,\zeta^{(0)}|=|Im\,\zeta^{(0)}|=1$, and $\zeta^{(1)}=\mathcal{O}(h)$ as $h\rightarrow 0$. Now we will construct two such $\zeta$'s explicitly. Let $\xi,\mu^{(1)},\mu^{(2)}\in\mathbb{R}^{n}$ be three real vectors with $|\mu^{(1)}|=|\mu^{(2)}|=1$ and $\xi\cdot\mu^{(1)}=\xi\cdot\mu^{(2)}=\mu^{(1)}\cdot\mu^{(2)}=0$. Set
\begin{equation}\label{zeta}
\zeta_{1}:=\displaystyle\frac{ih\xi}{2}+i\sqrt{1-h^{2}\frac{|\xi|^{2}}{4}}\mu^{(1)}+\mu^{(2)}, \quad \zeta_{2}:=-\displaystyle\frac{ih\xi}{2}+i\sqrt{1-h^{2}\frac{|\xi|^{2}}{4}}\mu^{(1)}-\mu^{(2)}.
\end{equation}
It is easy to check that $\zeta_{1},\zeta_{2}$ satisfy the above conditions with $\zeta^{(0)}_{1}=i\mu^{(1)}+\mu^{(2)}$ and $\zeta^{(0)}_{2}=i\mu^{(1)}-\mu^{(2)}$.

First we construct $u_{1}\in \mathcal{W}_{l_{2}}(\Sigma\cap B)$. To satisfy the boundary condition $u_{1}|_{l_{2}}=0$, we will reflect $\Sigma\cap B$ with respect to the boundary hyperplane $\Gamma_{2}$. More precisely, denote $(\Sigma\cap B)^{\ast}:=\{(x',-x_{n}):x=(x',x_{n})\in\Sigma\cap B\}$ where $x'=(x_{1},\cdots,x_{n-1})$. We extend the coefficients $A^{(1)}$ and $q^{(1)}$ to $(\Sigma\cap B)^{\ast}$ as follows: for the components $A^{(1)}_{j}$, $1\leq j\leq n-1$ and $q^{(1)}$, we extend them as even functions with respect to $\Gamma_{2}$, so we define
$$
\begin{array}{rl} \vspace{1ex}
\tilde{A}^{(1)}_{j}(x)=&
\left\{
\begin{array}{ll} \vspace{1ex}
A^{(1)}_{j}(x',x_{n}) & 0<x_{n}<L \\
A^{(1)}_{j}(x',-x_{n}) & -L<x_{n}<0 \\
\end{array}
\right.
, \quad j=1,\cdots,n-1 \\ \vspace{1ex}

\tilde{q}^{(1)}(x)=&
\left\{
\begin{array}{ll} \vspace{1ex}
q^{(1)}(x',x_{n}) & 0<x_{n}<L \\
q^{(1)}(x',-x_{n}) & -L<x_{n}<0 \\
\end{array}
\right. .\\
\end{array}
$$
For $A^{(1)}_{n}$ we extend it as an odd function with respect to $\Gamma_{2}$, that is, we define
$$
\begin{array}{rl} \vspace{1ex}
\tilde{A}^{(1)}_{n}(x)=&
\left\{
\begin{array}{ll} \vspace{1ex}
A^{(1)}_{n}(x',x_{n}) & 0<x_{n}<L \\
-A^{(1)}_{n}(x',-x_{n}) & -L<x_{n}<0 \\
\end{array}
\right.. \\
\end{array}
$$

On $\Gamma_{2}$, $A^{(1)}_{n}|_{x_{n}=0}=A^{(1)}\cdot\nu=0$ by \eqref{normal}, then $\tilde{A}^{(1)}\in L^{\infty}((\Sigma\cap B)\cup(\Sigma\cap B)^{\ast})$ and $\tilde{q}^{(1)}\in L^{\infty}((\Sigma\cap B)\cup(\Sigma\cap B)^{\ast})$. By Proposition \ref{existence}, there exist CGO solutions of the form
$$\tilde{u}_{1}(x,\zeta_{1},h)=e^{x\cdot\zeta_{1}/h}(e^{\Phi^{\sharp}_{1}(x,i\mu^{(1)}+\mu^{(2)},h)}+r_{1}(x,\zeta_{1},h))$$
which satisfy the equation $(\mathcal{L}_{\tilde{A}^{(1)},\tilde{q}^{(1)}}-k^{2})\tilde{u}_{1}=0$ in the bounded region $(\Sigma\cap B)\cup(\Sigma\cap B)^{\ast}$ with $\Phi^{\sharp}_{1}\in C^{\infty}(\overline{(\Sigma\cap B)\cup(\Sigma\cap B)^{\ast}})$,
\begin{equation}\label{a1}
(i\mu^{(1)}+\mu^{(2)})\cdot D \Phi^{\sharp}_{1}+(i\mu^{(1)}+\mu^{(2)})\cdot (\widetilde{A^{(1)}})^{\sharp}=0 \quad \textrm{ in } (\Sigma\cap B)\cup(\Sigma\cap B)^{\ast},
\end{equation}
\begin{equation}\label{a1derivative}
\|\partial^{\alpha}e^{\Phi^{\sharp}_{1}}\|_{L^{\infty}((\Sigma\cap B)\cup(\Sigma\cap B)^{\ast})}\leq C_{\alpha}h^{-|\alpha|/3} \textrm{ for all }  \alpha \textrm{ with } |\alpha|\geq 0,
\end{equation}
$\Phi^{\sharp}_{1}(\cdot,i\mu^{(1)}+\mu^{(2)},h)\rightarrow\Phi_{1}(\cdot,i\mu^{(1)}+\mu^{(2)})$ in $L^{2}((\Sigma\cap B)\cup (\Sigma\cap B)^{\ast})$ as $h\rightarrow 0$; and
\begin{equation}\label{r1}
\|r_{1}\|_{H^{1}_{scl}((\Sigma\cap B)\cup(\Sigma\cap B)^{\ast})}=\mathcal{O}(h^{1/3}) \textrm{ as } h\rightarrow 0.
\end{equation}
Let
\begin{equation} \label{CGOu1}
\begin{array}{rl}
u_{1}(x):= & \tilde{u}_{1}(x',x_{n})-\tilde{u}_{1}(x',-x_{n}) \\
 =& e^{x\cdot\zeta_{1}/h}(e^{\Phi^{\sharp}_{1}}+r_{1}(x))-\\
 & e^{(x',-x_{n})\cdot\zeta_{1}/h}(e^{\Phi^{\sharp}_{1}(x',-x_{n})}+r_{1}(x',-x_{n})) \quad\quad x\in\Sigma\cap B.\\
\end{array}
\end{equation}
It can be checked by direct computation that $u_{1}\in \mathcal{W}_{l_{2}}(\Sigma\cap B)$.\\

Next we construct $u_{2}\in \mathcal{V}_{l_{1}}(\Sigma\cap B)$. The construction is similar to that of $u_{1}$, but this time we reflect with respect to $\Gamma_{1}$. Denote $(\Sigma\cap B)^{\ast\ast}:=\{(x',-x_{n}+2L):x=(x',x_{n})\in\Sigma\cap B\}$ where $x'=(x_{1},\cdots,x_{n-1})$. We also extend the coefficients $A^{(2)}$ and $q^{(2)}$ to $(\Sigma\cap B)^{\ast\ast}$ as follows: for the components $A^{(2)}_{j},1\leq j\leq n-1$ and $q^{(2)}$, we extend them as even functions with respect to $\Gamma_{1}$, 
$$
\begin{array}{rl} \vspace{1ex}
\tilde{A}^{(2)}_{j}(x)=&
\left\{
\begin{array}{ll} \vspace{1ex}
A^{(2)}_{j}(x',x_{n}) & 0<x_{n}<L \\
A^{(2)}_{j}(x',-x_{n}+2L) & L<x_{n}<2L \\
\end{array}
\right.
, \quad j=1,\cdots,n-1 \\ \vspace{1ex}

\tilde{q}^{(2)}(x)=&
\left\{
\begin{array}{ll} \vspace{1ex}
q^{(2)}(x',x_{n}) & 0<x_{n}<L \\
q^{(2)}(x',-x_{n}+2L) & L<x_{n}<2L \\
\end{array}
\right. .\\
\end{array}
$$
For $A^{(2)}_{n}$ we extend it as an odd function with respect to $x_{n}=L$, i.e. we set
$$
\begin{array}{rl} \vspace{1ex}
\tilde{A}^{(2)}_{n}(x)=&
\left\{
\begin{array}{ll} \vspace{1ex}
A^{(2)}_{n}(x',x_{n}) & 0<x_{n}<L \\
-A^{(2)}_{n}(x',-x_{n}+2L) & L<x_{n}<2L \\
\end{array}
\right.. \\
\end{array}
$$

On $\Gamma_{1}$, $A^{(2)}_{n}|_{x_{n}=L}=A^{(2)}\cdot\nu=0$ by \eqref{normal}, then $\tilde{A}^{(2)}\in L^{\infty}((\Sigma\cap B)\cup(\Sigma\cap B)^{\ast\ast})$ and $\tilde{q}^{(2)}\in L^{\infty}((\Sigma\cap B)\cup(\Sigma\cap B)^{\ast\ast})$. By Proposition \ref{existence}, there exist CGO solutions of the form
$$\tilde{u}_{2}(x,\zeta_{2},h)=e^{x\cdot\zeta_{2}/h}(e^{\Phi^{\sharp}_{2}(x,i\mu^{(1)}-\mu^{(2)},h)}+r_{2}(x,\zeta_{2},h))$$
which satisfy the equation $(\mathcal{L}_{\overline{\tilde{A}^{(2)}},\overline{\tilde{q}^{(2)}}}-k^{2})\tilde{u}_{2}=0$ in the bounded region $(\Sigma\cap B)\cup(\Sigma\cap B)^{\ast\ast}$ with $\Phi^{\sharp}_{2}\in C^{\infty}(\overline{(\Sigma\cap B)\cup(\Sigma\cap B)^{\ast\ast}})$,
\begin{equation}\label{a2}
(i\mu^{(1)}-\mu^{(2)})\cdot D \Phi^{\sharp}_{2}+(i\mu^{(1)}-\mu^{(2)})\cdot (\widetilde{A^{(2)}})^{\sharp}=0 \quad \textrm{ in } (\Sigma\cap B)\cup(\Sigma\cap B)^{\ast\ast},
\end{equation}
\begin{equation}\label{a2derivative}
\|\partial^{\alpha}e^{\Phi^{\sharp}_{2}}\|_{L^{\infty}((\Sigma\cap B)\cup(\Sigma\cap B)^{\ast\ast})}\leq C_{\alpha}h^{-|\alpha|/3} \textrm{ for all }  \alpha \textrm{ with } |\alpha|\geq 0,
\end{equation}
$\Phi^{\sharp}_{2}(\cdot,i\mu^{(1)}-\mu^{(2)},h)\rightarrow\Phi_{2}(\cdot,i\mu^{(1)}-\mu^{(2)})$ in $L^{2}((\Sigma\cap B)\cup (\Sigma\cap B)^{\ast\ast})$ as $h\rightarrow 0$; and
\begin{equation}\label{r2}
\|r_{2}\|_{H^{1}_{scl}((\Sigma\cap B)\cup(\Sigma\cap B)^{\ast\ast})}=\mathcal{O}(h^{1/3}) \textrm{ as } h\rightarrow 0.
\end{equation}
Let
\begin{equation} \label{CGOu2}
\begin{array}{rl}
u_{2}(x):= & \tilde{u}_{2}(x',x_{n})-\tilde{u}_{2}(x',-x_{n}+2L) \\
 =& e^{x\cdot\zeta_{2}/h}(e^{\Phi^{\sharp}_{2}}+r_{2}(x))-\\
 & e^{(x',-x_{n}+2L)\cdot\zeta_{2}/h}(e^{\Phi^{\sharp}_{2}(x',-x_{n}+2L)}+r_{2}(x',-x_{n}+2L)) \quad x\in\Sigma\cap B.\\
\end{array}
\end{equation}
Then $u_{2}\in \mathcal{V}_{l_{1}}(\Sigma\cap B)$, which can be easily verified.\\

\subsection{END OF THE PROOF}
We are now in the position to finish the proof of Theorem 1.1. We will insert the partial data CGO solutions \eqref{CGOu1} and \eqref{CGOu2} into the identity \eqref{identity3}. For this purpose we compute
\begin{equation}
\begin{array}{rl} \vspace{1ex}
e^{x\cdot\zeta_{1}/h}e^{x\cdot\overline{\zeta}_{2}/h}=&e^{ix\cdot\xi} \\ \vspace{1ex}
e^{(x',-x_{n})\cdot\zeta_{1}/h}e^{x\cdot\overline{\zeta}_{2}/h}=&e^{-2\mu^{(2)}_{n}x_{n}/h+ib_{1}} \\ \vspace{1ex}
e^{x\cdot\zeta_{1}/h}e^{(x',-x_{n}+2L)\cdot\overline{\zeta}_{2}/h}=&e^{2\mu^{(2)}_{n}(x_{n}-L)/h+ib_{2}} \\ \vspace{1ex}
e^{(x',-x_{n})\cdot\zeta_{1}/h}e^{(x',-x_{n}+2L)\cdot\overline{\zeta}_{2}/h}=&e^{-2L\mu^{(2)}_{n}/h+ib_{3}} \\
\end{array}
\end{equation}
where $b_{1}, b_{2}, b_{3}\in\mathbb{R}^{n}$ are defined by
$$
\begin{array}{rl} \vspace{1ex}
b_{1}:=&x'\cdot\xi'-\displaystyle\frac{2}{h}\sqrt{1-h^{2}\frac{|\xi|^{2}}{4}}\mu^{(1)}_{n}x_{n}, \\ \vspace{1ex}
b_{2}:=&x'\cdot\xi'+\displaystyle\frac{2}{h}\sqrt{1-h^{2}\frac{|\xi|^{2}}{4}}\mu^{(1)}_{n}(x_{n}-L)+L\xi_{n}, \\ \vspace{1ex}
b_{3}:=&x'\cdot\xi'-\displaystyle\frac{2L}{h}\sqrt{1-h^{2}\frac{|\xi|^{2}}{4}}\mu^{(1)}_{n}-x_{n}\xi_{n}+L\xi_{n}. \\
\end{array}
$$
In order to eliminate the undesired terms, we shall further assume that $\mu^{(2)}_{n}>0$ so that for $0<x_{n}<L$ we have the following pointwise convergence:
\begin{equation}\label{vanish}
\begin{array}{rl} \vspace{1ex}
e^{(x',-x_{n})\cdot\zeta_{1}/h}e^{x\cdot\overline{\zeta}_{2}/h}&\rightarrow 0 \textrm{ as } h\rightarrow 0^+,\\ \vspace{1ex}
e^{x\cdot\zeta_{1}/h}e^{(x',-x_{n}+2L)\cdot\overline{\zeta}_{2}/h}&\rightarrow 0 \textrm{ as } h\rightarrow 0^+,\\ \vspace{1ex}
e^{(x',-x_{n})\cdot\zeta_{1}/h}e^{(x',-x_{n}+2L)\cdot\overline{\zeta}_{2}/h}&\rightarrow 0 \textrm{ as } h\rightarrow 0^+.\\
\end{array}
\end{equation}
Notice that by \eqref{a1derivative}, \eqref{r1}, \eqref{a2derivative} and \eqref{r2}, we have for $j=1,2,$
\begin{equation}\label{order}
\begin{array}{cc}
\|e^{\Phi^{\sharp}_{j}}\|_{L^{\infty}(\Sigma\cap B)}=\mathcal{O}(1) & \|De^{\Phi^{\sharp}_{j}}\|_{L^{\infty}(\Sigma\cap B)}=\mathcal{O}(h^{-1/3}) \\
\|r_{j}\|_{L^{2}(\Sigma\cap B)}=\mathcal{O}(h^{1/3}) & \|Dr_{j}\|_{L^{2}(\Sigma\cap B)}=\mathcal{O}(h^{-2/3}). \\
\end{array}
\end{equation}
Therefore, with the complex geometric optics solutions $u_{1}$ and $u_{2}$ given by \eqref{CGOu1} and \eqref{CGOu2}, we conclude from \eqref{vanish} and \eqref{order} that, after multiplying the identity \eqref{identity3} by $h$, the second term will tend to zero as $h\rightarrow 0^+$, i.e.
\begin{equation}\label{term2}
h\displaystyle\int_{\Sigma\cap B}((A^{(1)})^{2}-(A^{(2)})^{2}+q^{(1)}-q^{(2)})u_{1}\overline{u_{2}}\,dx\rightarrow 0 \quad \textrm{ as } h\rightarrow 0^+.
\end{equation}

Now we analyze the first term after multiplying \eqref{identity3} by $h$. Denote $\zeta^{\ast}_{j}=(\zeta'_{j},-(\zeta_{j})_{n})$ for $\zeta_{j}=(\zeta'_{j},(\zeta_{j})_{n})$, $j=1,2$. Using \eqref{CGOu1} and \eqref{CGOu2} we compute
\begin{equation}
\begin{array}{rl} \vspace{1ex}
Du_{1}(x)=&-\displaystyle\frac{i\zeta_{1}}{h}e^{x\cdot\zeta_{1}/h}(e^{\Phi^{\sharp}_{1}(x)}+r_{1}(x))+e^{x\cdot\zeta_{1}/h}(D e^{\Phi^{\sharp}_{1}(x)}+Dr_{1}(x))\\ \vspace{1ex}
&+\displaystyle\frac{i\zeta^{\ast}_{1}}{h}e^{(x',-x_{n})\cdot\zeta_{1}/h}(e^{\Phi^{\sharp}_{1}(x',-x_{n})}+r_{1}(x',-x_{n})) \\ \vspace{1ex}
&-e^{(x',-x_{n})\cdot\zeta_{1}/h}(D e^{\Phi^{\sharp}_{1}(x',-x_{n})}+Dr_{1}(x',-x_{n})).\\
\end{array}
\end{equation}
\begin{equation}
\begin{array}{rl} \vspace{1ex}
\overline{Du_{2}}(x)=&\displaystyle\frac{i\overline{\zeta_{2}}}{h}e^{x\cdot\overline{\zeta_{2}}/h}(e^{\overline{\Phi^{\sharp}_{2}(x)}}+\overline{r_{2}(x)})+e^{x\cdot\overline{\zeta_{2}}/h}(\overline{D e^{\Phi^{\sharp}_{2}(x)}}+\overline{Dr_{2}(x)})\\ \vspace{1ex}
&-\displaystyle\frac{i\overline{\zeta^{\ast}_{2}}}{h}e^{(x',-x_{n}+2L)\cdot\overline{\zeta_{2}}/h}(e^{\overline{\Phi^{\sharp}_{2}(x',-x_{n}+2L)}}+\overline{r_{2}(x',-x_{n}+2L)}) \\ \vspace{1ex}
&-e^{(x',-x_{n}+2L)\cdot\overline{\zeta_{2}}/h}(\overline{D e^{\Phi^{\sharp}_{2}(x',-x_{n}+2L)}}+\overline{Dr_{2}(x',-x_{n}+2L)}).\\
\end{array}
\end{equation}
Combining the above computation with the facts that $\Phi^{\sharp}_{1}(\cdot,i\mu^{(1)}+\mu^{(2)},h)\rightarrow \Phi_{1}(\cdot,i\mu^{(1)}+\mu^{(2)})$ in $L^{2}((\Sigma\cap B)\cup(\Sigma\cap B)^{\ast})$ as $h\rightarrow 0^+$ and $\Phi^{\sharp}_{2}(\cdot,i\mu^{(1)}-\mu^{(2)},h)\rightarrow \Phi_{2}(\cdot,i\mu^{(1)}-\mu^{(2)})$ in $L^{2}((\Sigma\cap B)\cup(\Sigma\cap B)^{\ast\ast})$ as $h\rightarrow 0^+$, we have that, as $h\rightarrow 0^+$,
\begin{equation}\label{term1}
\begin{array}{rl}\vspace{1ex}
 &h\displaystyle\int_{\Sigma\cap B}(A^{(1)}-A^{(2)})\cdot Du_{1}\overline{u_{2}}\,dx \\ \vspace{1ex}
 \rightarrow& (i\mu^{(1)}+\mu^{(2)})\cdot\displaystyle\int_{\Sigma\cap B}(A^{(1)}-A^{(2)}) e^{ix\cdot\xi}e^{\Phi_{1}(x,i\mu^{(1)}+\mu^{(2)})+\overline{\Phi_{2}(x,i\mu^{(1)}-\mu^{(2)})}}\,dx.\\
\end{array}
\end{equation}
Therefore, multiplying \eqref{identity3} by $h$ and letting $h\rightarrow 0^+$ for the constructed solutions $u_{1}$ and $u_{2}$, we obtain from \eqref{term2} and \eqref{term1} that
$$
(i\mu^{(1)}+\mu^{(2)})\cdot\displaystyle\int_{\Sigma\cap B}(A^{(1)}-A^{(2)}) e^{ix\cdot\xi}e^{\Phi_{1}(x,i\mu^{(1)}+\mu^{(2)})+\overline{\Phi_{2}(x,i\mu^{(1)}-\mu^{(2)})}}\,dx=0.
$$
In fact, \cite[Lemma 6.2]{S2} implies that the same identity is true with $e^{\Phi_{1}(x,i\mu^{(1)}+\mu^{(2)})+\overline{\Phi_{2}(x,i\mu^{(1)}-\mu^{(2)})}}$ replaced by $1$, i.e.
\begin{equation}
(i\mu^{(1)}+\mu^{(2)})\cdot\displaystyle\int_{\Sigma\cap B}(A^{(1)}-A^{(2)}) e^{ix\cdot\xi}\,dx=0
\end{equation}
for all $\xi,\mu^{(1)},\mu^{(2)}\in\mathbb{R}^{n}$ satisfying
$$\xi\cdot\mu^{(1)}=\xi\cdot\mu^{(2)}=\mu^{(1)}\cdot\mu^{(2)}=0,\quad |\mu^{(1)}|=|\mu^{(2)}|=1,\quad \mu^{(2)}_{n}>0.$$

This implies the vanishing of the Fourier transform of components of the distribution $d(A^{(1)}-A^{(2)})$ in the first quadrant $\{\xi\in\mathbb{R}^{n}:\xi_{1}>0,\cdots,\xi_{n}>0\}$, and hence in $\mathbb{R}^{n}$ by analyticity of the Fourier transform of compactly supported distributions. For details of this type of argument we refer to \cite{KLU}.

Therefore, $d(A^{(1)}-A^{(2)})=0$ in $\Sigma$. Since $\Sigma$ is simply connected, there exists compactly supported $\Psi\in W^{1,\infty}(\overline{\Sigma})$ such that
$$A^{(1)}-A^{(2)}=\nabla\Psi \quad \textrm{ in } \Sigma.$$
Moreover, it can be shown as in \cite{KLU} that $\Psi=0$ on $\partial\Sigma$. Thus, we can replace $A^{(1)}$ by $A^{(1)}+\nabla\Psi$ while keeping the Dirichlet-to-Neumann map unchanged due to \eqref{gauge}. In the following we will assume this replacement has been made so that $A^{(1)}=A^{(2)}$.

Inserting $A^{(1)}=A^{(2)}$ into \eqref{identity3} yields
$$\displaystyle\int_{\Sigma\cap B}(q^{(1)}-q^{(2)})u_{1}\overline{u_{2}}\,dx=0.$$
Again we plug in the CGO solutions $u_{1}$ as in \eqref{CGOu1} and $u_{2}$ as in \eqref{CGOu2}, then take $h\rightarrow 0^+$ to obtain
$$\displaystyle\int_{\Sigma\cap B}(q^{(1)}-q^{(2)})e^{ix\cdot\xi}e^{\Phi_{1}(x)+\overline{\Phi_{2}(x)}}\,dx=0.$$
As before, by using the argument in \cite[Lemma 6.2]{S2}, the above identity is still true after replacing $e^{\Phi_{1}(x)+\overline{\Phi_{2}(x)}}$ by $1$, which gives
$$\displaystyle\int_{\Sigma\cap B}(q^{(1)}-q^{(2)})e^{ix\cdot\xi}\,dx=0.$$
This identity is valid for all $\xi$ for which there exist $\mu^{(1)},\mu^{(2)}\in\mathbb{R}^{n}$ with
$$\xi\cdot\mu^{(1)}=\xi\cdot\mu^{(2)}=\mu^{(1)}\cdot\mu^{(2)}=0,\quad |\mu^{(1)}|=|\mu^{(2)}|=1,\quad \mu^{(2)}_{n}>0.$$
It is clear that for $\xi$ in the first quadrant, the vectors $\mu^{(1)},\mu^{(2)}$ always exist. This implies that the Fourier transform of $(q^{(1)}-q^{(2)})\chi_{\Sigma\cap B}$ vanishes in the first quadrant, hence in $\mathbb{R}^{n}$ by the analyticity of the Fourier transform. Here $\chi_{\Sigma\cap B}$ denotes the characteristic function of the set $\Sigma\cap B$. This completes the proof of Theorem 1.1.

\section{Proof of Theorem 1.2}
The proof of Theorem 1.2 is analogous to that of Theorem 1.1. We will derive an integral identity, construct some partial data CGO solutions, and finally derive the vanishing of the Fourier transform of some compactly supported distributions to conclude the uniqueness.

First, by a similar argument as in Section 3.1, we can obtain the identity \eqref{identity3} for all $u_{1}\in \mathcal{W}_{l_{2}}(\Sigma\cap B)$ and $u_{2}\in \mathcal{V}_{l_{2}}(\Sigma\cap B)$. We will continue to use complex frequencies $\zeta_{1}$ and $\zeta_{2}$ defined in \eqref{zeta} as well as $u_{1}$ of the form \eqref{CGOu1}.

To construct $u_{2}\in \mathcal{V}_{l_{2}}(\Sigma\cap B)$, we proceed as in the definition of $u_{1}$ by reflecting the coefficients with respect to the boundary hyperplane $\Gamma_{2}$. For $A^{(2)}_{j}$, $1\leq j\leq n-1$ and $q^{(2)}$, we extend them as even functions with respect to $x_{n}=0$ and define 
$$
\begin{array}{rl} \vspace{1ex}
\tilde{A}^{(2)}_{j}(x)=&
\left\{
\begin{array}{ll} \vspace{1ex}
A^{(2)}_{j}(x',x_{n}) & 0<x_{n}<L \\
A^{(2)}_{j}(x',-x_{n}) & -L<x_{n}<0 \\
\end{array}
\right.
, \quad j=1,\cdots,n-1 \\ \vspace{1ex}

\tilde{q}^{(2)}(x)=&
\left\{
\begin{array}{ll} \vspace{1ex}
q^{(2)}(x',x_{n}) & 0<x_{n}<L \\
q^{(2)}(x',-x_{n}) & -L<x_{n}<0 \\
\end{array}
\right. .\\
\end{array}
$$

For $A^{(2)}_{n}$, we extend it as an odd function with respect to $\Gamma_{2}$ and define
$$
\begin{array}{rl} \vspace{1ex}
\tilde{A}^{(2)}_{n}(x)=&
\left\{
\begin{array}{ll} \vspace{1ex}
A^{(2)}_{n}(x',x_{n}) & 0<x_{n}<L \\
-A^{(2)}_{n}(x',-x_{n}) & -L<x_{n}<0 \\
\end{array}
\right.. \\
\end{array}
$$

On $\Gamma_{2}$, we have $A^{(2)}_{n}|_{x_{n}=0}=A^{(2)}\cdot\nu=0$ by \eqref{normal}, then $\tilde{A}^{(2)}\in L^{\infty}((\Sigma\cap B)\cup(\Sigma\cap B)^{\ast})$ and $\tilde{q}^{(2)}\in L^{\infty}((\Sigma\cap B)\cup(\Sigma\cap B)^{\ast})$. Proposition \ref{existence} implies that there exist CGO solutions of the form
$$\tilde{u}_{2}(x,\zeta_{2},h)=e^{x\cdot\zeta_{2}/h}(e^{\Phi^{\sharp}_{2}(x,i\mu^{(1)}-\mu^{(2)},h)}+r_{2}(x,\zeta_{2},h))$$
which satisfy the equation $(\mathcal{L}_{\overline{\tilde{A}^{(2)}},\overline{\tilde{q}^{(2)}}}-k^{2})u_{2}=0$ in the bounded region $(\Sigma\cap B)\cup(\Sigma\cap B)^{\ast}$ with $\Phi^{\sharp}_{2}\in C^{\infty}(\overline{(\Sigma\cap B)\cup(\Sigma\cap B)^{\ast}})$,
\begin{equation}\label{a3}
(i\mu^{(1)}-\mu^{(2)})\cdot D \Phi^{\sharp}_{2}+(i\mu^{(1)}-\mu^{(2)})\cdot (A^{(2)})^{\sharp}=0 \quad \textrm{ in } (\Sigma\cap B)\cup(\Sigma\cap B)^{\ast},
\end{equation}
\begin{equation}\label{a3derivative}
\|\partial^{\alpha}e^{\Phi^{\sharp}_{2}}\|_{L^{\infty}((\Sigma\cap B)\cup(\Sigma\cap B)^{\ast})}\leq C_{\alpha}h^{-|\alpha|/3} \textrm{ for all }  \alpha \textrm{ with } |\alpha|\geq 0,
\end{equation}
$\Phi^{\sharp}_{2}(\cdot,i\mu^{(1)}-\mu^{(2)},h)\rightarrow\Phi_{2}(\cdot,i\mu^{(1)}-\mu^{(2)})$ in $L^{2}((\Sigma\cap B)\cup (\Sigma\cap B)^{\ast})$ as $h\rightarrow 0$; and
\begin{equation}\label{r3}
\|r_{2}\|_{H^{1}_{scl}((\Sigma\cap B)\cup(\Sigma\cap B)^{\ast})}=\mathcal{O}(h^{1/3}) \textrm{ as } h\rightarrow 0.
\end{equation}
Let
\begin{equation} \label{CGOu3}
\begin{array}{rl}
u_{2}(x):= & \tilde{u}_{2}(x',x_{n})-\tilde{u}_{2}(x',-x_{n}) \\
= & e^{x\cdot\zeta_{2}/h}(e^{\Phi^{\sharp}_{2}(x)}+r_{2}(x))- \\
 & e^{(x',-x_{n})\cdot\zeta_{2}}(e^{\Phi^{\sharp}_{2}(x',-x_{n})}+r_{2}(x',-x_{n})) \quad\quad x\in\Sigma\cap B.
\end{array}
\end{equation}
Then $u_{2}\in \mathcal{V}_{l_{2}}(\Sigma\cap B)$.\\

We will insert $u_{1}$ in \eqref{CGOu1} and $u_{2}$ in \eqref{CGOu3} into the identity \eqref{identity3}. To this end, we compute some products which appear in the integral
$$
\begin{array}{rl} \vspace{1ex}
e^{x\cdot\zeta_{1}/h}e^{x\cdot \overline{\zeta}_{2}/h}=&e^{ix\cdot\xi} \\ \vspace{1ex}
e^{x\cdot\zeta_{1}/h}e^{(x',-x_{n})\cdot\overline{\zeta}_{2}/h}=&e^{ix\cdot c_{1}+2\mu^{(2)}_{n}x_{n}/h}\\ \vspace{1ex}
e^{(x',-x_{n})\cdot\zeta_{1}/h}e^{x\cdot\overline{\zeta}_{2}/h}=&e^{ix\cdot c_{2}-2\mu^{(2)}_{n}x_{n}/h} \\ \vspace{1ex}
e^{(x',-x_{n})\cdot\zeta_{1}/h}e^{(x',-x_{n})\cdot\overline{\zeta}_{2}/h}=&e^{i(x',-x_{n})\cdot\xi} \\
\end{array}
$$
where
$$
c_{1}=\left(\xi',+\frac{2}{h}\sqrt{1-h^{2}\frac{|\xi|^{2}}{4}}\mu^{(1)}_{n}\right), \quad c_{2}=\left(\xi',-\frac{2}{h}\sqrt{1-h^{2}\frac{|\xi|^{2}}{4}}\mu^{(1)}_{n}\right).
$$
To eliminate the undesired terms, we assume $\mu^{(2)}_{n}=0$ and $\mu^{(1)}_{n}\neq 0$, so $c_{1},c_{2}\rightarrow\infty$ as $h\rightarrow 0$. We have
$$
\begin{array}{rl}\vspace{1ex}
 & \zeta_{1}\cdot\displaystyle\int_{\Sigma\cap B} (A^{(1)}-A^{(2)})e^{x\cdot\zeta_{1}/h}e^{(x',-x_{n})\cdot\overline{\zeta_{2}}/h}e^{\Phi^{\sharp}_{1}(x)+
\overline{\Phi^{\sharp}_{2}(x',-x_{n})}}\,dx \\ \vspace{1ex}
=& \zeta_{1}\cdot\displaystyle\int_{\Sigma\cap B} (A^{(1)}-A^{(2)})e^{ix\cdot\tilde{\xi}_{+}}e^{\Phi_{1}(x)+
\overline{\Phi_{2}(x',-x_{n})}}\,dx \\ \vspace{1ex}
 &+\zeta_{1}\cdot\displaystyle\int_{\Sigma\cap B} (A^{(1)}-A^{(2)})e^{ix\cdot\tilde{\xi}_{+}}(e^{\Phi^{\sharp}_{1}(x)+
\overline{\Phi^{\sharp}_{2}(x',-x_{n})}}-e^{\Phi_{1}(x)+
\overline{\Phi_{2}(x',-x_{n})}})\,dx \rightarrow 0 \\
\end{array}
$$
as $h\rightarrow 0$. Here the first integral on the right hand side tends to zero by the Riemann-Lebesgue lemma; the second tends to zero since $\Phi^{\sharp}_{1}(\cdot,i\mu^{(1)}+\mu^{(2)},h)\rightarrow \Phi_{1}(\cdot,i\mu^{(1)}+\mu^{(2)})$ in $L^{2}((\Sigma\cap B)\cup(\Sigma\cap B)^{\ast})$ as $h\rightarrow 0^+$, $\Phi^{\sharp}_{2}(\cdot,i\mu^{(1)}-\mu^{(2)},h)\rightarrow \Phi_{2}(\cdot,i\mu^{(1)}-\mu^{(2)})$ in $L^{2}((\Sigma\cap B)\cup(\Sigma\cap B)^{\ast})$ as $h\rightarrow 0^+$, and since the inequality
$$|e^{z}-e^{w}|\leq |z-w|e^{\textrm{max}(Re\,z,Re\,w)}, \quad\quad z,w\in\mathbb{C}^{n}.$$
Similarly
\begin{equation}
\zeta_{1}\cdot\displaystyle\int_{\Sigma\cap B} (A^{(1)}-A^{(2)})e^{(x',-x_{n})\cdot \zeta_{1}/h}e^{x\cdot\overline{\zeta}_{2}/h}e^{\Phi^{\sharp}_{1}(x)+
\overline{\Phi^{\sharp}_{2}(x',-x_{n})}}\,dx\rightarrow 0
\end{equation}
as $h\rightarrow 0$. Therefore, multiplying \eqref{identity3} by $h$ and letting $h\rightarrow 0$ we get
$$
\left(i\mu^{(1)}+\mu^{(2)}\right)\cdot\displaystyle\int_{(\Sigma\cap B)\cup(\Sigma\cap B)^{\ast}}\left(\tilde{A}^{(1)}-\tilde{A}^{(2)}\right)e^{ix\cdot\xi}e^{\Phi_{1}(x)+\overline{\Phi_{2}(x)}}\,dx=0
$$
where we have made a change of variable so that the integral domain becomes $(\Sigma\cap B)\cup(\Sigma\cap B)^{\ast}$. As before we may replace $e^{\Phi_{1}(x)+\overline{\Phi_{2}(x)}}$ by $1$ to obtain
$$
\left(i\mu^{(1)}+\mu^{(2)}\right)\cdot\displaystyle\int_{(\Sigma\cap B)\cup(\Sigma\cap B)^{\ast}}\left(\tilde{A}^{(1)}-\tilde{A}^{(2)}\right)e^{ix\cdot\xi}\,dx=0
$$
for all $\xi,\mu^{(1)},\mu^{(2)}\in\mathbb{R}^{n}$ such that
$$
\xi\cdot\mu^{(1)}=\xi\cdot\mu^{(2)}=\mu^{(1)}\cdot\mu^{(2)}=0,\quad |\mu^{(1)}|=|\mu^{(2)}|=1, \quad \mu^{(2)}_{n}=0, \quad \mu^{(1)}_{n}\neq 0.
$$
This implies the vanishing of the Fourier transform of components of the distribution $d(\tilde{A}^{(1)}-\tilde{A}^{(2)})$ in the first quadrant, and hence in $\mathbb{R}^{n}$ by analyticity of the Fourier transform of compactly supported distributions. We conclude that $d\tilde{A}^{(1)}=d\tilde{A}^{(2)}$.

The rest of the proof is similar to that of Theorem 1.1. $d\tilde{A}^{(1)}=d\tilde{A}^{(2)}$ implies the existence of a function $\Psi\in W^{1,\infty}({\overline{\Sigma}})$ with $\Psi=0$ along $\partial((\Sigma\cap B)\cup(\Sigma\cap B)^{\ast})$ such that $\tilde{A}^{(1)}-\tilde{A}^{(2)}=\nabla\Psi$. Replacing $\tilde{A}^{(1)}$ by $\tilde{A}^{(1)}+\nabla\Psi$ if necessary, we may assume that $\tilde{A}^{(1)}=\tilde{A}^{(2)}$. Inserting this into
\eqref{identity3} and arguing as in the proof of Theorem 1.1 we arrive at
$$\displaystyle\int_{(\Sigma\cap B)\cup(\Sigma\cap B)^{\ast}}(\tilde{q}^{(1)}-\tilde{q}^{(2)})e^{ix\cdot\xi}\,dx=0$$
for all $\xi$ such that there exist $\mu^{(1)},\mu^{(2)}\in\mathbb{R}^{n}$ with
$$
\xi\cdot\mu^{(1)}=\xi\cdot\mu^{(2)}=\mu^{(1)}\cdot\mu^{(2)}=0,\quad |\mu^{(1)}|=|\mu^{(2)}|=1, \quad \mu^{(2)}_{n}=0, \quad \mu^{(1)}_{n}\neq 0.
$$
Since for any $\xi$ in the first quadrant we can find such vectors $\mu^{(1)}$ and $\mu^{(2)}$, we conclude that the Fourier transform of $(\tilde{q}^{(1)}-\tilde{q}^{(2)})\chi_{(\Sigma\cap B)\cup(\Sigma\cap B)^{\ast}}$ vanishes
in the first quadrant, hence in $\mathbb{R}^{n}$ by analyticity. This completes the proof of Theorem 1.2

%
%
%
%

\end{document}